\documentclass{article}
\begin{document}
\centerline{\textbf{\Large On Lehmer Binomial Series}}
\[
\]
\centerline{\bf Nikos Bagis}
\centerline{nikosbagis@hotmail.gr}
\[
\]
\centerline{\textbf{Abstract}}
We give evaluations in closed form of certain Lehmer-binomial series 
\[
\]
\section{Lehmer Series and roots of polynomials}

We begin with the observation that the integral 
\begin{equation}
\int^{1}_{0}\frac{dt}{1-(1-t)^{s-1}tz}
\end{equation}
have expansion
\begin{equation}
\sum_{\rho}\frac{\log(1-\rho)-\log(-\rho)}{z\rho^{s-2}\left(1-s+s\rho\right)}
\end{equation}
where $\sum_{\rho/H}g(\rho)$ denotes the sum of $g$ over all roots of the equation $H(X)=0$. Here we have $H(X)=1-zX^{s-1}+zX^s$ hence the equation is
$$
1-zX^{s-1}+zX^s=0\eqno{:(eq)} 
$$
Since we get all the roots of $(eq)$ in the sum (2), we don't have to worry about denoting them with some order (if the roots are simple).\\
The integral (1) is equal to the Lehmer's hypergeometric series
\begin{equation}
A_s(z):=\sum^{\infty}_{m=1}\frac{z^m}{(sm+1)\left(^{sm}_m\right)}
\end{equation}
This can be proved by expanding the binomial in (3) into $\Gamma$ functions and use $$B(a,b)=\int^{1}_{0}t^{a-1}(1-t)^{b-1}dt=\frac{\Gamma(a)\Gamma(b)}{\Gamma(a+b)}.$$
After suming the series we get the desired result. By this way we lead to the following\\
\\
\textbf{Theorem 1.}\\
Set $p_{s,k}(z)$ to be the $k$-root (in any descrete order) of the equation (1), then
\begin{equation}
A_s(z)=\sum^{\infty}_{m=1}\frac{z^m}{(sm+1)\left(^{sm}_m\right)}= 
\frac{1}{z}\sum^{s}_{k=1}\frac{\log\left(1-p_{s,k}(z)\right)-\log\left(-p_{s,k}(z)\right)}{p_{s,k}(z)^{s-2}\left(1-s+sp_{s,k}(z)\right)}
\end{equation}
\textbf{Proof.}\\
For to prove the equality of (1) and (2) we first observe that the differential equation
\begin{equation}
y'(x)=\sum^{\infty}_{k=1}a_ky(x)^k=H(y(x))
\end{equation}
have solution
\begin{equation}
x+C=\sum_{\rho/H}\frac{\log(y(x)-\rho)}{H'(\rho)} ,
\end{equation}
where the sumation is taken over all roots of $H(x)=0$ (here $H$ is a simple roots polynomial). If we invert $y$ we get 
$$
y^{(-1)}(x)=\int\frac{1}{H(x)}dx=\sum_{\rho/H}\frac{\log(x-\rho)}{H'(\rho)}
$$  
But with the change of variable $t\rightarrow (1-t)$ relation (1) becomes
$$
\int^{1}_{0}\frac{dt}{1-t^{s-1}(1-t)z}=\int^{1}_{0}\frac{dt}{1-zt^{s-1}+zt^s}
$$
Hence the implication is proved, since $H'(x)=zx^{s-2}(1-s+sx)$.\\
\\
\textbf{Examples}\\
\textbf{1.} For $s=2$ the two roots of $(eq)$ are $X_{1,2}=\frac{1}{2}\pm\frac{1}{2}\sqrt{\frac{-4+z}{z}}$. Setting $z\rightarrow e^w$ and derivating two times the Lehmer series with respect to $w$ we get a large function. Setting $w=\log(2)$ and simplifing we get
\begin{equation}
\sum^{\infty}_{n=1}\frac{n^22^n}{(2n+1)\left(^{2n}_n\right)}=  1+\frac{\pi}{2}
\end{equation}
\textbf{2.}
$$
12\sum^{\infty}_{n=1}\frac{n^22^{-n}}{(2n+1)\left(^{2n}_n\right)}=
$$
\begin{equation}
={}_3F_{2}\left[2,2,2;1,5/2;1/8\right]=\frac{12}{343} \left(-56+26 \sqrt{7} \pi -52 \sqrt{7} \arctan\left[\frac{3}{\sqrt{7}}\right]\right)
\end{equation}
\textbf{3.}
$$
\sum^{\infty}_{n=0}\frac{(-1)^{n}}{(3n+1)4^n\left(^{3n}_n\right)}=
$$
\begin{equation}
={}_3F_2\left[1/2,1,1;2/3,4/3;-1/27\right]=-\frac{5 \textrm{arccot}\left(\frac{5}{\sqrt{7}}\right)}{2\sqrt{7}}+\frac{3\log(2)}{4}
\end{equation}
\\
In the same way as above we can prove\\

\section{A general case of Lehmer series and the Polylogarithm function}

\textbf{Theorem 2.}
\begin{equation}
\sum^{\infty}_{n=0}\frac{n^kx^n}{(ns+1)\left(^{ns}_n\right)}=\int^{1}_{0}L_i\left(-k,x(1-t)t^{s-1}\right)dt
\end{equation}
where $L_i(k,z)=\sum^{\infty}_{n=1}\frac{z^n}{n^k}$ is the polylogarithm function.\\
\\

Note that the integral in (10) can evaluated in a closed form with Mathematica. A closed type form as in Theorem 1 is pointless and quite complicated. Some evaluations directly with Mathematica are\\
\\
\textbf{4.}
\begin{equation}
\sum^{\infty}_{n=0}\frac{n^3}{(3n+1)2^n\left(^{3n}_n\right)}=\frac{1335-22 \pi -351 \log(2)}{15625}
\end{equation} 
\textbf{5.}
\begin{equation}
\sum^{\infty}_{n=0}\frac{n}{(2n+1)\left(^{2n}_n\right)}=-\frac{2}{27} \left(-9+\sqrt{3} \pi \right)
\end{equation}
\textbf{6.} 
Let $H(X)=3-X^2+X^3$, then
\begin{equation}
\frac{1}{36} {}_5F_{4}\left[\left\{\frac{3}{2},2,2,2,2\right\},\left\{1,1,\frac{5}{3},\frac{7}{3}\right\},\frac{4}{81}\right]=\frac{30033}{456533}+
$$
$$
+\frac{3}{456533}\sum_{\rho/H}\frac{-13897\left[\log(1-\rho)-\log(-\rho)\right]+2738\rho\left[\log(1-\rho)-\log(-\rho)\right]}{-2\rho+3\rho^2}
\end{equation}

\section{References}

F.J. Dyson, Norman E Frankel, M.L. Glasser. 'Lehmer's Interesting Series'. Amer. Math. Monthly. Vol. 120,No. 2,pp.116-130. (2013)

\end{document}